\documentclass{aims}
\usepackage{amsmath}
  \usepackage{paralist}
  \usepackage{graphics} 
  \usepackage{epsfig} 
\usepackage{graphicx}  \usepackage{epstopdf}
 \usepackage[colorlinks=true]{hyperref}
\hypersetup{urlcolor=blue, citecolor=red}

\def\currentvolume{X}
 \def\currentissue{X}
  \def\currentyear{200X}
   \def\currentmonth{XX}
    \def\ppages{X--XX}
     \def\DOI{10.3934/xx.xx.xx.xx}

\newtheorem{theorem}{Theorem}[section]

\theoremstyle{definition}

\newtheorem{remark}{Remark}

\def\bk{{\bf{k}}}

\def\bj{{\bf{j}}}

\title[] 
      {Spectral Convergence of the Stochastic Galerkin Approximation to the Boltzmann Equation with Multiple Scales
and Large Random Perturbation in the Collision Kernel
}

\author[Esther S. Daus, Shi Jin and Liu Liu]{}


 \email{esther.daus@tuwien.ac.at}
 \email{shijin-m@sjtu.edu.cn}
 \email{lliu@ices.utexas.edu}

\thanks{The first author acknowledges partial support from the
 Austrian Science Fund (FWF), grants P27352 and P30000, the second author is supported by NSF grants DMS-1522184, DMS-1819012 and DMS-1107291: RNMS KI-Net, NSFC grant No.31571071 and No.11871297, the third author is supported by the funding
DOE--Simulation Center for Runaway Electron Avoidance and Mitigation, project number DE-SC0016283. }

\begin{document}
\maketitle

\centerline{\scshape Esther S. Daus$^*$}
\medskip
{\footnotesize
 \centerline{Institute for Analysis and Scientific Computing, Vienna University of Technology, }
  \centerline{Wiedner Hauptstrasse 8--10, 1040 Wien, Austria}
} 
\medskip

\centerline{\scshape Shi Jin$^*$}
\medskip
{\footnotesize
 \centerline{School of Mathematical Sciences, Institute of Natural Sciences, MOE-LSC and SHL-MAC, }
  \centerline{Shanghai Jiao Tong University, Shanghai 200240, China}
}

\medskip
\centerline{\scshape Liu Liu$^*$}
\medskip
{\footnotesize
   \centerline{Institute for Computational Engineering and Sciences, }
   \centerline{University of Texas at Austin, Austin, Texas 78705, USA}
}

\bigskip


\begin{abstract}
In [L. Liu and S. Jin, {\it Multiscale Model. Simult.}, 16, 1085--1114, 2018],
    spectral convergence and long-time decay of the numerical solution towards the global equilibrium of the stochastic Galerkin approximation for the
    Boltzmann equation with random inputs in the initial data and
    collision kernel for hard potentials and Maxwellian molecules under Grad's angular cutoff were established using the hypocoercive properties of the collisional   kinetic model. One
    assumption for the random perturbation of the collision kernel
    is that the perturbation is in the order of the
    Knudsen number, which can be very small in the fluid dynamical regime.
    In this article, we remove this smallness assumption, and establish
    the same results but now for random perturbations of the collision kernel
    that can be of order one. The new analysis relies on the establishment of a spectral gap for the numerical collision operator. 
\end{abstract}

\section{Introduction}

Kinetic equations are usually derived from $n$-body Newton's equations
via the mean-field limit \cite{BGP}. As such, inevitably they contain
uncertainties in the initial or boundary data, forcing or source terms, and
in particular, in their collision kernels or scattering cross-sections.
Quantifying such uncertainties are important to assess, validate and improve
the kinetic models. For some recent efforts in uncertainty quantification
for kinetic equations, see for example \cite{HuReview, Twoband, QinWang, JinLiuMa, Liu18}.  

In \cite{LiuJin-UQ}, by extending the hypocoercivity theory developed for
general deterministic  collisional nonlinear kinetic equations \cite{MB, CN} to the random uncertain setting, the authors established the regularity and
long-time behavior of the solution with random initial data and, for
the case of the Boltzmann equation, random collision kernel as well, under
suitable assumptions.
Moreover, for the stochastic Galerkin (SG) approximation of the random Boltzmann
equation, the spectral convergence
and long-time exponential error decay were also established, under the condition
that the random perturbation of the collision kernel is in the order of
the Knudsen number. In the fluid dynamical regime, this assumption becomes
quite restrictive since the Knudsen number can be vanishingly small.
In this paper, we establish the same results by removing this constraint,
namely we allow the
random perturbation of the collision kernel to be of {\it order one}.

The improved result is based on the observation that the SG
system of the random Boltzmann equation bears some similarity
with the multi-species Boltzmann equation. In \cite{HY-MB}, the hypocoercivity framework for the linearized multi-species Boltzmann equation was established.
One essential ingredient in the analysis is to establish the spectral-gap estimate of the linearized multi-species collision operator. The argument that
leads to such a spectral-gap estimate, using the symmetries of the Boltzmann collision
operator, can be modified for the linearized
collision operator for the SG system, which, despite of the anisotropic
nature of the (numerical) collision kernel, can also have a spectral
gap thanks to a careful vectorial handling of the collision operator for the approximated functions using ideas similar to \cite{HY-MB}. 
In the previous analysis \cite{LiuJin-UQ}, the off-diagonal entries of the
(numerical) collision kernel, which consist of the randomly perturbed
part of the collision kernel, were absorbed into the nonlinear part,
creating a perturbation of order $O(\epsilon)$ of the linear part of the collision
terms. In this work, the restriction that the random perturbation needs to be of order $O(\epsilon)$ is removed.

In Section \ref{sec:Intro}, we introduce the (uncertain) Boltmann equation and its
basic properties.
Section \ref{sec:gPC} presents the SG method for solving the Boltzmann
 equations with random  uncertainties, and reviews the previous results in \cite{LiuJin-UQ}.
 In Section \ref{sec:main}, we establish the main results, namely the proof
of the spectral gap for the linearized numerical collision operator, and we explain
how the $\mathcal O(\epsilon)$ random perturbation assumption for the collision kernel can be removed.

\section{Introduction of the Boltzmann equation with uncertainties}
\label{sec:Intro}

We first give a review of some of the results in \cite{LiuJin-UQ} that will be useful in our proof.
Consider the initial value problem for the Boltzmann equation
\begin{align}
\label{model}
\left\{
\begin{array}{l}
\displaystyle \partial_t f + \frac{1}{\epsilon^\alpha} v\cdot\nabla_x f  =
\frac{1}{\epsilon^{1+\alpha}} \mathcal Q(f, f), \\[4pt]
\displaystyle  f(0,x,v,z)= f_{I}(x,v,z), \qquad 
x\in\mathbb T^{d_x},  \, v\in \mathbb R^{d_v}, 
z\in I_z\subset \mathbb R, 
\end{array}\right.
\end{align}
where $f=f(t,x,v,z)$ is the particle density distribution  that depends on time $t$,
particle position $x\in\mathbb T^{d_x}$ (periodic box of $d_x$ dimension), velocity $v\in\mathbb R^{d_v}$ and a random variable $z$.
The numbers $d_x, d_v\geq 1$ denote the dimension of the spatial and velocity spaces, and $z$ is a random variable that lies in the domain 
$I_z\subset \mathbb R$ with compact support, which is used to account for the random uncertainties or inputs.
The operator $\mathcal Q$ is quadratic and models the binary collisional interactions between particles. The parameter $\epsilon$ is the dimensionless Knudsen number, the 
ratio of particle mean free path over the domain size. The choice $\alpha=1$ refers to the incompressible Navier-Stokes scaling, 
and $\alpha=0$ corresponds
to the Euler (or acoustic in this article) scaling.
Moreover, we assume periodic boundary conditions on the torus $\mathbb T^{d_x}$. 

For notational simplicity, we set $d_v=3$ in the following. Since we consider random collision kernels, the operator $\mathcal Q$ is defined by
$$\mathcal Q(f, g) =\int_{\mathbb R^3 \times \mathbb S^2} B(|v-v^{\ast}|, \cos\theta, z)
(f^{\prime}g^{\prime\ast} - f g^{\ast}) dv^{\ast}d\sigma, $$
where we used the abbreviations $f^{\prime}=f(v^{\prime})$, $g^{\ast}=g(v^{\ast})$ and $g^{\prime\ast}=g(v^{\prime\ast})$, and
$\mathbb S^2$ is the three-dimensional unit sphere. Note that
$v^{\prime}$ and $v^{\prime\ast}$ are the post-collisional velocities of particles depending on the pre-collisional velocities
$v$ and $v^{\ast}$.  During elastic collisions, the momentum and
kinetic energy of the involved particles are conserved, namely,
$$ v^{\prime}=\frac{v+v^{\ast}}{2}+\frac{|v-v^{\ast}|}{2}\sigma, \qquad
v^{\prime\ast}=\frac{v+v^{\ast}}{2}-\frac{|v-v^{\ast}|}{2}\sigma, $$
where 
$\sigma \in \mathbb{S}^2$ is a parameter on the 2-dimensional unit sphere.
The collision kernel $B=B(|v-v^{\ast}|, \cos\theta, z)$
 is a non-negative function depending on the modulus of the relative velocity $|v-v^{\ast}|$, the cosinus of the deviation angle $\theta$ with
$$\cos\theta=\sigma\cdot(v-v^{\ast})/|v-v^{\ast}|, $$ and the random variable $z\in I_z$. 
\\[8pt]
\indent {\bf Properties of the collision operators:}
First, conservation of mass, momentum and energy is satisfied, \textit{i.e.}
$$ \frac{d}{dt} \int_{\mathbb T^{d_x}\times\mathbb R^3}\,  f \begin{pmatrix} 1 \\ v \\ |v|^2 \end{pmatrix} dvdx =
\frac{d}{dt}\int_{\mathbb T^{d_x}\times\mathbb R^3}\,  \mathcal Q(f, f) \begin{pmatrix} 1 \\ v \\ |v|^2 \end{pmatrix} dvdx = 0.  $$
Next, we have the dissipation of entropy
\begin{equation} \int_{\mathbb R^3}\, \mathcal Q(f, f)\ln f\, dv \leq 0,
\end{equation}
which is known as the celebrated Boltzmann's H-theorem. Moreover,
\begin{equation} \int_{\mathbb R^3}\, \mathcal Q(f, f)\ln f\, dv=0  \, \Leftrightarrow\,
\mathcal Q(f, f)=0  \, \Leftrightarrow\, f = \mathcal M_{\text{loc}}, \end{equation}
where $\mathcal M_{\text{loc}}$ is the {\it local} equilibrium state given by a Maxwellian distribution
$$\mathcal M_{\text{loc}} =  \frac{\rho}{(2\pi T)^{3/2}}\exp\left(-\frac{|v-u|^2}{2T}\right), $$ 
with 
\begin{equation}\label{Rho} \rho =\int_{\mathbb R^3} f\, dxdv, \qquad u =\frac{1}{\rho}\int_{\mathbb R^3} f v\, dxdv, \qquad
\, T =\frac{1}{3\rho}\int_{\mathbb R^3} f |v-u|^2\, dxdv. 
\end{equation}

The {\it global} equilibrium is the unique stationary solution to (\ref{model}) and is given by
\begin{equation}\label{M}\mathcal M(v) =\frac{1}{(2\pi)^{\frac{3}{2}}}\, e^{-\frac{|v|^2}{2}}, \end{equation} 
where by translating and scaling the coordinate system, we assumed $\rho=1$, $u=0$ and $T=1$ in (\ref{M}).
For further properties of the Boltzmann equation, see \cite{Cer}.

One of the central questions in kinetic theory is to understand the long-time
behavior of the solution, and for this, the hypocoercive effects of the kinetic equations play a
pivotal role, see \cite{CD, Guo-NS, HeNi, Strain-Guo1, Villani}. A hypocoercivity
framework for generic nonlinear collisional kinetic equations was
established in \cite{CN, MB}. In \cite{LiuJin-UQ}, this framework was extended to nonlinear collisional kinetic equations with random initial data and/or random collision
kernels, which allow to study the long-time sensitivity, regularity,
and exponential decay of the (random) solution towards the (deterministic)
global equilibrium, for both the random kinetic equations and their stochastic
Galerkin approximations. Note that these studies have been carried out only for solutions
near the global equilibrium, \textit{i.e.} in a {\it perturbative} setting, such that the solution can be defined under suitable Sobolev norms. Let
\begin{equation} \label{per-f}
f =\mathcal M + \epsilon\sqrt{\mathcal M}\, h. 
\end{equation}
Now inserting this ansatz into the model (\ref{model}), the fluctuation $h$ satisfies
\begin{align}
&\displaystyle
\label{h-eq}
\left\{
\begin{array}{l}
\displaystyle  \partial_t h + \frac{1}{\epsilon^{\alpha}}v\cdot\nabla_x h = \frac{1}{\epsilon^{1+\alpha}}\mathcal L(h)
+ \frac{1}{\epsilon^{\alpha}}\mathcal F(h, h),   \\[4pt]
 \displaystyle h(t=0) = h_I,
 \end{array}\right.
\end{align}
where the linearized collision operator is defined by
\begin{align}
&\displaystyle\mathcal L(h)= \mathcal M^{-1/2}\left[\mathcal Q(\sqrt{\mathcal M}\, h,\mathcal M) + 
\mathcal Q(\mathcal M, \sqrt{\mathcal M}\, h)\right] \notag\\[4pt]
&\displaystyle \qquad = \mathcal M^{1/2} \int_{\mathbb R^3\times\mathbb S^2}B(|v-v^*|, \cos\theta,z)\mathcal M^{\ast} \notag\\[4pt]
&\displaystyle\label{L0}\qquad\qquad\qquad \cdot\left[\frac{h^{\ast\prime}}{(\mathcal M^{\ast\prime})^{1/2}}+\frac{h^{\prime}}{(\mathcal M^{\prime})^{1/2}}\
-\frac{h^{\ast}}{(\mathcal M^{\ast})^{1/2}}-\frac{h}{\mathcal M^{1/2}}\right] dv^{\ast}d\sigma, 
\end{align}
while the nonlinear operator has the form
\begin{align}
  &\displaystyle\mathcal F(h,h)=\mathcal M^{-1/2}\left[\mathcal Q(\sqrt{\mathcal M}\, h, \sqrt{\mathcal M}\, h) +
  \mathcal Q(\sqrt{\mathcal M}\, h, \sqrt{\mathcal M}\, h)\right] \notag\\[2pt]
  &\displaystyle\label{GB}\qquad\quad =\int_{\mathbb R^3\times\mathbb S^2}
    B(|v-v^*|, \cos\theta,z)(\mathcal M^{\ast})^{1/2}(h^{\ast\prime}h^{\prime}-h^{\ast}h)\, dv^{\ast}d\sigma. 
\end{align}

The linearized operator $\mathcal L$ is acting on $L^2_v=\{f\, |\int_{\mathbb R^3} f^2\, dv <\infty\}$,
with a finite dimensional kernel $N(\mathcal L)=\text{Span}\{\varphi_1, \cdots, \varphi_n\}$, where
$\{\varphi_i\}_{1\leq i\leq n}$ is an orthonormal family of polynomials in $v$ corresponding to the manifold of local equilibria for the linearized kinetic models.
The orthogonal projection on $N(\mathcal L)$ in $L^2_v$ is defined by
\begin{equation}\label{Pi}\Pi_{\mathcal L} (h) = \sum_{i=1}^{n}\, \left(\int_{\mathbb R^3} h\varphi_i\, dv\right)\varphi_i. \end{equation}
In the classical case of hard potentials and Maxwellian molecules under Grad's angular cutoff, meaning that the collision kernel $B$ has the form 
$$B(|v-v^*|, \cos\theta,z)=\Phi(| v-v^*|)b(\cos\theta,z)$$ 
with the kinetic part $\Phi$ satisfying  
\begin{align}\label{kinetic.part}
\Phi(|v-v^*|)=C|v-v^*|^{\gamma} ~~ for ~~\gamma \in [0,1]
\end{align}
and the angular part $b$ being locally integrable, assumptions (1.1), (1.2), (1.3) in \cite{MC} are satisfied. Thus,  \cite[Theorem 1.1]{MC} holds, and consequently $\mathcal L$ has the following {\it local coercivity property with an explicitly computable coercivity constant $\lambda>0$}: 
There exists a constant $\lambda>0$, such that $\forall\, h\in L^2_v$,
\begin{equation}\label{coercivity} \langle\mathcal L(h),\, h\rangle_{L^2_v} \leq -\lambda\, ||h^{\perp}||_{\Lambda_v}^2, \end{equation}
where
\begin{equation} h^{\perp}=h-\Pi_{\mathcal L}(h)
\end{equation} stands for the microscopic part of $h$, and the coercivity norm is
$$||h||_{\Lambda_v}=||h(1+|v|)^{\gamma/2}||_{L^2_v}, $$
where 
$\gamma$ denotes the constant $\gamma\in [0,1]$ in \eqref{kinetic.part}. For all $z$, 
we make the same assumption for $b(\cos\theta, z)$ as in \cite[assumption (1.3)]{MC}, namely: 
$$\inf_{\sigma_1, \sigma_2\in\mathbb S^2}\int_{\sigma_3\in\mathbb S^2}\min\{b(\sigma_1\cdot\sigma_3, z), b(\sigma_2\cdot\sigma_3, z)\}
d\sigma_3 >0. $$

\section{A gPC based Stochastic Galerkin Method}
\label{sec:gPC}

We first review the gPC-SG method for solving kinetic equations with uncertainties \cite{Hu}. For more general differential equations, see for examples \cite{Ghanem, Xiu}. We want to approximate $f$ (or $h$) in the following way
\begin{align}
&\displaystyle  f(t,x,v,z)\approx \sum_{|\bk|=1}^{K}\, f_{\bk}(t,x,v)\psi_{\bk}(z) := f^K(t,x,v,z), \notag\\[2pt]
&\displaystyle \label{gPC-Ans} h(t,x,v,z)\approx \sum_{|\bk|=1}^{K}\, h_{\bk}(t,x,v)\psi_{\bk}(z) := h^K(t,x,v,z).
\end{align}
Here $\bk=(k_1, \cdots, k_n)$ is a multi-index with $|\bk|=k_1+\cdots + k_n$.
The orthonormal gPC basis functions $\{\psi_{\bk}(z)\}$ satisfy
$$\int_{I_z}\, \psi_{\bk}(z)\psi_{\bj}(z)\pi(z)dz=\delta_{\bk\bj}, \qquad 1\leq |\bk|, \, |\bj|\leq K, $$
where $\pi(z)$ is the probability distribution function of $z$, which is given
{\it a priori} in our problem. 
Note that $f$ can be expanded by
$$f(t,x,v,z)=\sum_{|\bk|=1}^{\infty} \hat f_{\bk}(t,x,v)\psi_{\bk}(z), \quad \textnormal{where}~~
\hat f_{\bk}(t,x,v)=\int_{I_z}\, f(t,x,v,z)\psi_{\bk}(z)\pi(z)dz. $$
Finally, we define the projection operator $P_{K}$ as
\begin{equation}\label{proj} P_{K} f(t,x,v,z)=\sum_{|\bk|=1}^K\, \hat f_{\bk}(t,x,v)\psi_{\bk}(z). \end{equation}

We will consider the one-dimensional random variable $z$ in the sequel.
By inserting the ansatz (\ref{gPC-Ans}) into (\ref{h-eq})  and conducting a standard Galerkin projection,
one obtains the following gPC-SG system for $h_k$ for each $1\leq k \leq K$:
\begin{align}
\label{h-gPC}
\left\{
\begin{array}{l}
\displaystyle \displaystyle \partial_t h_k + \frac{1}{\epsilon^\alpha}v\cdot\nabla_x h_k =\frac{1}{\epsilon^{1+\alpha}}\mathcal L_k(h^K) +
\displaystyle \frac{1}{\epsilon^\alpha}\mathcal F_k(h^K, h^K),  \\[10pt]
\displaystyle  h_k(0,x,v)=h_k^{0}(x,v),
\end{array}\right.
\end{align}
with periodic boundary conditions and the initial data
$$h_k^{0} := \int_{I_z} h^{0}(x,v,z)\psi_k(z)\pi(z)dz. $$
In (\ref{h-gPC}), the operator $\mathcal L_k(h^K)$ is given by
\begin{align}
&\displaystyle \mathcal L_k(h^K) = \langle \mathcal L(h^K), \, \psi_k \rangle_{L^2(\pi(z))} = \langle \mathcal L(\sum_{j=1}^K h_j \psi_j), \, \psi_k \rangle_{L^2(\pi(z))} \notag\\[4pt]
&\displaystyle\qquad\quad = \sum_{j=1}^K\, \mathcal M^{1/2} \int_{\mathbb R^3 \times \mathbb S^2} S_{kj}\, 
\mathcal M^{\ast}\notag\\[4pt]
&\label{L00}\displaystyle\qquad\qquad\qquad\qquad\quad 
\cdot\left(\frac{h_j^{\prime\ast}}{(\mathcal M^{\prime\ast})^{1/2}} 
+ \frac{h_j^{\prime}}{(\mathcal M^{\prime})^{1/2}} - \frac{h_j^{\ast}}{(\mathcal M^{\ast})^{1/2}} - \frac{h_j}{\mathcal M^{1/2}}\right) dv^{\ast}d\sigma,
\end{align}
where we denote ${\bf S}_{K\times K}$ by the $K\times K$ matrix with $(k, j)$-th component
\begin{equation}\label{SS} S_{kj} = \int_{I_z} B\, \psi_k \psi_j \pi(z)dz = \int_{I_z} B\, \psi_k \psi_j\, d\pi(z), \qquad 1\leq k, j \leq K. \end{equation}
Note that $S$ is symmetric. 

The nonlinear term is
\begin{align*}
&\displaystyle \mathcal F_k(h^K, h^K) = \langle \mathcal F(h^K, h^K), \, \psi_k \rangle_{L^2(\pi(z))}\\[4pt]
&\displaystyle \qquad\qquad\quad = \sum_{i, j=1}^{K}\int_{\mathbb R^3\times\mathbb S^2} T_{kij}\, (\mathcal M^{\ast})^{1/2}
\left(h_i^{\ast}\, h_j^{\ast\prime} - h_i\, h_j^{\ast}\right) dv^{\ast}d\sigma,
\end{align*}
with the tensor ${\bf T}_{K\times K\times K}$ defined by
$$T_{kij}=\int_{I_z} B\, \psi_k(z)\psi_i(z)\psi_j(z)\pi(z)dz. $$
\\[4pt]
\indent We first review \cite[Theorem 5.1]{LiuJin-UQ}, which is the main result on the SG systems of that paper. 
\begin{theorem}
\label{thm:5.1}
Assume that the collision kernel $B$ satisfies the assumptions
\begin{align}
&\displaystyle B(|v-v^{\ast}|, \cos\theta, z) = \Phi(|v-v^{\ast}|)b(\cos\theta,z), \quad
\phi(\xi) = C_{\phi}\, \xi^{\gamma} \text{ with }\gamma\in[0,1], C_{\phi}>0, \notag\\[6pt]
&\label{BK1}\displaystyle \forall\eta\in[-1,1], \qquad |b(\eta,z)|\leq C_b, \, |\partial_{\eta} b(\eta,z)|\leq \tilde C_b, \,
|\partial_z^k b(\eta, z)|\leq C_b^{\ast},  \,  \forall\, 0\leq k\leq r,
\end{align}
where $b$ is linear in $z$, and has the particular form
\begin{equation}\label{B-LN} b(\cos\theta, z)= b_0(\cos\theta)+  b_1(\cos\theta)z,  \end{equation}
where $z\in I_z$ has a compact support, that is, 
\begin{equation}\label{z-BD} |z| \leq C_z, \end{equation} and
\begin{equation}\label{Assump-B} |\partial_z b| = |b_1| \leq \mathcal O(\epsilon).
\end{equation}
We also assume the technical condition \cite{ShuJin}
\begin{equation}\label{basis}
||\psi_k||_{L^{\infty}} \leq C k^p, \qquad 1\leq k\leq K,  \end{equation}
with a parameter $p \geq 0$. Let $q>p+2$, and define the energy $E^{K}$ by
\begin{equation}\label{Def-E} E^{K}(t) = E_{s,q}^{K}(t) = \sum_{k=1}^{K}\, ||k^q h_k||_{H_{x,v}^s}^2, \end{equation}
with the initial data satisfying
\begin{equation}\label{E-IC} E^{K}(0) \leq \eta. \end{equation}
Finally, let for all $s\geq s_0$, $0\leq \epsilon_d\leq 1$ and
for $0\leq \epsilon\leq \epsilon_d$, let $h^K$ be a gPC solution of (\ref{h-gPC}) with
\begin{equation}\label{h-IC} ||h_{\text{in}}||_{H_{x,v}^{s,r}L_z^{\infty}} \leq C_I. \end{equation}
Then we get that
\begin{align*}
&\displaystyle E^{K}(t)\leq \eta\, e^{-\varepsilon^{1-\alpha}\,\tau t}\,,  \\[6pt]
&\displaystyle ||h^K||_{H^s_{x,v}L_z^{\infty}} \leq \tilde\eta\, e^{-\varepsilon^{1-\alpha}\,\tau t}, \qquad ||h^K||_{H^s_{x,v}L_z^2} \leq \tilde\eta\,
e^{-\varepsilon^{1-\alpha}\,\tau t}\,,  \\[6pt]
&\displaystyle ||h-h^K||_{H_{x,v}^{s} L_z^2} \leq C\, \frac{e^{-\varepsilon^{1-\alpha}\, \tilde\tau t}}{K^r}\,,
\end{align*}
where $\eta$, $\tilde\eta$, $\tau$, $\tilde\tau$ and $C$ are all positive constants independent of $K$ and $z$.
\end{theorem}
Note that the theorems on the estimate of $E^K$, the gPC error and spectral
convergence of the SG method in \cite{LiuJin-UQ} require the same set of assumptions on $B$ and gPC polynomial basis $\psi$,
given by (\ref{BK1})--(\ref{basis}). We point out that the linearity hypothesis (\ref{B-LN}) on $b$ resembles the form of the well known Karhunen-Loeve expansion, which is widely used for modeling random fields and is linear in $z$. 
One may also put the random dependence on $\phi$ in the collision kernel $B$, which brings a similar analysis and the same conclusion.

\section{The main result}
\label{sec:main}

The main purpose of this section is to obtain the same results as in Theorem \ref{thm:5.1}, whereas
relaxing the condition (\ref{Assump-B}) from $|\partial_z b| \leq \mathcal O(\varepsilon)$ to $|\partial_z b| = \mathcal O(1)$.

We will focus on the linearized operator $\mathcal L_k$ below.
We denote \begin{equation}
\label{LB}\mathcal L^{B}(h) := \mathcal M^{1/2} \int_{\mathbb R^3 \times \mathbb S^2} B \mathcal M^{\ast}\left(\frac{h^{\prime\ast}}{(\mathcal M^{\prime\ast})^{1/2}}
+ \frac{h^{\prime}}{(\mathcal M^{\prime})^{1/2}} - \frac{h^{\ast}}{(\mathcal M^{\ast})^{1/2}} - \frac{h}{\mathcal M^{1/2}}\right) dv^{\ast}d\sigma,
\end{equation}
and we introduce the operator $\Theta$ by
\begin{equation}\label{theta} \Theta[\tilde h] = \tilde h^{\prime\ast} + \tilde h^{\prime} - \tilde h^{\ast} - \tilde h, \end{equation}
with $\tilde h = h/\mathcal M^{1/2}$.
Then (\ref{LB}) can be rewritten as
\begin{equation}\label{LB1}\mathcal L^{B}(h) := \mathcal M^{1/2} \int_{\mathbb R^3 \times \mathbb S^2} B \mathcal M^{\ast}\, \Theta[\tilde h]\, dv^{\ast}d\sigma.
\end{equation}

To remove the $\mathcal O(\varepsilon)$ assumption on the collision kernels, we need to reinvestigate the estimate for the term involving the linearized operator
$\mathcal L_k$, namely
\begin{equation}\label{LL} \sum_{k=1}^K\, k^{2q}\, \langle \mathcal L_k(h^K), \, h_k \rangle_{L^2_v}.  \end{equation}

By (\ref{L0}), 
\begin{equation}\label{LS} \mathcal L_k(h^K) =
 \sum_{j=1}^K\, \mathcal M^{1/2} \int_{\mathbb R^3 \times \mathbb S^2} \underbrace{\int_{I_z} B\, \psi_k \psi_j \pi(z)dz}_{S_{kj}}
\mathcal M^{\ast}\, \Theta[\tilde h_j]\, dv^{\ast}d\sigma =  \sum_{j=1}^K\, \mathcal L^{S_{kj}}(h_j), \end{equation}
where $$  \mathcal L^{S_{kj}}(h_j)=\mathcal M^{1/2} \int_{\mathbb R^3 \times \mathbb S^2} S_{kj}\mathcal M^{\ast}\, \Theta[\tilde h_j]\, dv^{\ast}d\sigma, $$
by the definition (\ref{LB1}). 

\begin{remark}
In \cite{LiuJin-UQ}, we assume that $B$ satisfies (\ref{BK1})--(\ref{Assump-B}), that is,
$$ B(|v-v^{\ast}|, \cos\theta, z)=\phi(|v-v^{\ast}|) b(\cos\theta, z) = \phi(|v-v^{\ast}|)\left(b_0(\cos\theta) + b_1(\cos\theta)z\right), $$
with $|\partial_z b| = |b_1| \leq \mathcal O(\varepsilon)$. Thus $|S_{kk}| \leq C_b$, using (\ref{z-BD}), for $k \neq j$ we have
\begin{align*}
&\displaystyle  |S_{kj}| = \left| b_1 \int_{I_z} z\, \psi_k \psi_j \pi(z)dz\right| \leq |b_1|\, ||z||_{L^{\infty}} \langle |\psi_k|,  |\psi_j| \rangle_{L^2(\pi(z))} \\[4pt]
&\displaystyle  \qquad\leq |b_1|\, C_z\, ||\psi_k||_{L^2(\pi(z))} ||\psi_j||_{L^2(\pi(z))} = |b_1|\, C_z \leq \mathcal O(\varepsilon).
\end{align*}
The procedure below makes it possible to remove the assumption that the random perturbation in the collision kernels require a bound of order
$\mathcal O(\varepsilon)$, and avoids using the assumption that $|S_{kj}| \leq C \varepsilon$ for $k\neq j$ as it was shown in \cite{LiuJin-UQ}.
\end{remark}

Consider the term
\begin{align}
&\displaystyle \text{Term I}: = \sum_k\, k^{2q}\, \langle \mathcal L_k(h^K), \, h_k \rangle_{L^2_v}
= \sum_k \sum_j\, k^{2q}\, \langle \mathcal L^{S_{kj}}(h_j), \, h_k \rangle_{L^2_v} \notag\\[4pt]
&\displaystyle \qquad\quad = \sum_k \sum_j\,  k^{2q}\int_{\mathbb R^3} \int_{\mathbb R^3 \times \mathbb S^2} \mathcal M^{1/2}\, S_{kj}\, \mathcal M^{\ast}\, \Theta[\tilde h_j] h_k\, dv^{\ast}d\sigma dv \notag\\[4pt]
&\displaystyle\label{E0} \qquad\quad = \sum_k \sum_j\, k^{2q}\int_{\mathbb R^3} \int_{\mathbb R^3 \times \mathbb S^2} S_{kj}\, \mathcal M \mathcal M^{\ast}\underbrace{\mathcal M^{-1/2} h_k}_{\tilde h_k}\, \Theta[\tilde h_j]\, dv^{\ast} d\sigma dv.
\end{align}
Thus (for simplicity, we omit writing the integral domain below),
\begin{equation}\label{E1} \text{Term I} =  \sum_k \sum_j\, k^{2q}\int \int S_{kj}\, \mathcal M \mathcal M^{\ast}\, \Theta[\tilde h_j]\, \tilde h_k \, dv^{\ast} d\sigma dv.
\end{equation}
\\[2pt]
{\bf Step 1:} We make the change of variables $(v, v^{\ast},\sigma) \rightarrow (v^{\prime}, v^{\prime\ast},k)$ with $k=(v-v^*)/|v-v^*|$ on the right hand side of (\ref{E1}), which has unit Jacobian and is involutive. Note that $S_{kj}$ also depends on $|v-v^{\ast}|$ through $B$, with the relation $|v-v^{\ast}| = |v^{\prime} - v^{\prime\ast}|$, thus in the following steps $S_{kj}$ will not
be  affected by the coordinate transformation. Then one obtains
\begin{align}
&\displaystyle \text{Term I} =  \sum_k \sum_j\, k^{2q}\int \int S_{kj}\, \mathcal M^{\prime} \mathcal M^{\prime\ast}
\left( -\Theta[\tilde h_j]\right) \tilde h_k^{\prime}\, dv^{\ast} d\sigma dv \notag\\[4pt]
&\displaystyle \label{E2}
\qquad\quad = - \sum_k \sum_j\, k^{2q}\int \int S_{kj}\, \mathcal M \mathcal M^{\ast}\, \Theta[\tilde h_j]\, \tilde h_k^{\prime}\, dv^{\ast} d\sigma dv,
\end{align}
where we used $\mathcal M^{\prime} \mathcal M^{\prime\ast} = \mathcal M \mathcal M^{\ast}$, due to the conservation of kinetic energy, i.e.,
$|v|^2 + |v^{\ast}|^2 = |v^{\prime}|^2 + |v^{\prime\ast}|^2$.
\\[8pt]
{\bf Step 2:} Then one makes the change of variables $(v, v^{\ast}) \rightarrow (v^{\ast}, v)$ on the right hand side of (\ref{E1}) and obtains
\begin{equation}\label{E3}
\text{Term I} = \sum_k \sum_j\, k^{2q}\int \int S_{kj}\, \mathcal M \mathcal M^{\ast}\, \Theta[\tilde h_j]\, \tilde h_k^{\ast}\, dv^{\ast}d\sigma dv.
\end{equation}
\\[8pt]
{\bf Step 3:} Now we use again the change of variables $(v, v^{\ast},\sigma) \rightarrow (v^{\prime}, v^{\prime\ast},k)$ with $k=(v-v^*)/|v-v^*|$
on the right hand side of (\ref{E3}), thus
\begin{align}
&\displaystyle \text{Term I} = \sum_k \sum_j\, k^{2q}\int \int S_{kj}\, \mathcal M^{\prime} \mathcal M^{\prime\ast}
 \left( -\Theta[\tilde h_j] \right) \tilde h_k^{\prime\ast}\, dv^{\ast}d\sigma dv \notag\\[4pt]
&\displaystyle \label{E4} \qquad\quad = - \sum_k \sum_j\, k^{2q}\int \int S_{kj}\, \mathcal M \mathcal M^{\ast}\, \Theta[\tilde h_j]\, \tilde h_k^{\prime\ast}\,  dv^{\ast}d\sigma dv,
\end{align}
where $\mathcal M^{\prime} \mathcal M^{\prime\ast} = \mathcal M \mathcal M^{\ast}$ is used.
\\[8pt]
{\bf Step 4:} Finally, by combining (\ref{E1}), (\ref{E2}), (\ref{E3}) and (\ref{E4}), one has
\begin{align}
&\displaystyle \text{Term I} = \frac{1}{4} \sum_k \sum_j\, k^{2q}\int \int S_{kj}\, \mathcal M \mathcal M^{\ast}\, \Theta[\tilde h_j]
\left(\tilde h_k - \tilde h_k^{\prime} + \tilde h_k^{\ast} - \tilde h_k^{\prime\ast}\right) dv^{\ast}d\sigma dv \notag\\[4pt]
&\displaystyle \qquad\quad = \frac{1}{4}\sum_k \sum_j\, k^{2q}\int\int S_{kj}\, \mathcal M \mathcal M^{\ast}\, \Theta[\tilde h_j] \left(- \Theta[\tilde h_k]\right)dv^{\ast}d\sigma dv \notag\\[4pt]
&\displaystyle\label{E5} \qquad\quad  = - \frac{1}{4}\sum_k \sum_j\, k^{2q}\int\int S_{kj}\, \mathcal M \mathcal M^{\ast}\, \Theta[\tilde h_j]\, \Theta[\tilde h_k]\, dv^{\ast}d\sigma dv.
\end{align}
Recall the definition of operator $\Theta$ in (\ref{theta}) and matrix ${\bf S}$ in (\ref{SS}).
Denote the vector ${\bf \Theta}_{K\times 1}$ with $j$-th component given by
$$ \Theta_j := \tilde h_j^{\prime\ast} + \tilde h_j^{\prime} - \tilde h_j^{\ast} - \tilde h_j, \qquad 1\leq j \leq K. $$
Then (\ref{E5}) can be written as
\begin{align}
&\displaystyle\label{E6a} \text{Term I} = - \frac{1}{4}\sum_k \sum_j\, k^{2q}\int \int  \mathcal M \mathcal M^{\ast}\, S_{kj}\, \Theta_j\, \Theta_k\, dv^{\ast}d\sigma dv\\[4pt]
&\displaystyle
\qquad\quad = -\frac{1}{4}\int\int\mathcal M\mathcal M^{\ast}\sum_k \sum_j \left(\frac{k}{j}\right)^{q}\, S_{kj}\left(j^q\, \Theta_j\right)\left(k^q\, \Theta_k\right)dv^{\ast}d\sigma dv \notag\\[4pt]
&\displaystyle\label{E6b}
\qquad\quad = -\frac{1}{4}\int\int\mathcal M\mathcal M^{\ast}\sum_k \sum_j \left(\frac{k}{j}\right)^{q}\, S_{kj}\, \tilde\Theta_j\, \tilde\Theta_k\, dv^{\ast}d\sigma dv, 
\end{align}
where we denote $$\tilde\Theta_j=j^q\, \Theta, \qquad 1\leq j\leq K. $$

Now we focus on the summation
\begin{equation}\label{ST}\sum_k \sum_j \left(\frac{k}{j}\right)^{q}\, S_{kj}\, \tilde\Theta_j\, \tilde\Theta_k
\end{equation} in the integral (\ref{E6b}). 
Recall the definition of the matrix $S$ and the assumptions (\ref{BK1})--(\ref{B-LN}) on the collision kernel in Theorem \ref{thm:5.1}. Then we have that the $(k,j)$-th component of $S$ is given by 
\begin{align*}
&\displaystyle S_{kj}(|v-v^{\ast}|, \cos\theta) = \int_{I_z}\phi(|v-v^{\ast}|) b(\cos\theta, z)\psi_k \psi_j\, d\pi(z) \\[4pt]
&\displaystyle\qquad\qquad\qquad\qquad = \phi(|v-v^{\ast}|) \int_{I_z}(b_0(\cos\theta) + b_1(\cos\theta)z) \psi_k \psi_j\, d\pi(z). 
\end{align*}

We make the following assumption on the collision kernel in addition to (\ref{BK1}), (\ref{B-LN}) and (\ref{z-BD}). 
For a fixed $\theta\in[0, \pi]$, 
$\exists\, D(\cos\theta)>0$, such that
\begin{equation}\label{BK2} b_0 (\cos\theta) \geq (2^q+2)\, |b_1(\cos\theta)|\, C_z + D(\cos\theta), \end{equation}
where $D(\cos\theta)$ satisfies the same assumption as $b(\cos\theta)$ in \cite[assumption (1.3)]{MC}. 
Denote the matrix $\tilde S_{K\times K}$ by (we omit writing the $\cos\theta$ dependence in $\tilde S$, $b_0$ and $b_1$ below)
$$ \tilde S_{kj} = \int_{I_z}(b_0 + b_1 z)\, \psi_k \psi_j\, d\pi(z) = b_0\, \delta_{kj} + b_1 \int_{I_z} z\, \psi_k\psi_j\, d\pi(z). $$
Note that $\tilde S$ is a tridiagonal matrix, see $(3.4)$ in \cite{Xiu-Shen}, i.e., for a fixed $k$, $\tilde S_{kj}\neq 0$ only when 
either $j=k-1$, $j=k$ or $j=k+1$ happens. 
\\[2pt]

We now consider the following term: 
\begin{align*}
&\displaystyle \text{Term A}=\sum_k \sum_j \left(\frac{k}{j}\right)^{q}\, \tilde S_{kj}\, \tilde\Theta_j\, \tilde\Theta_k \\[4pt]
&\displaystyle\qquad\quad  = b_0 \sum_k\sum_j \left(\frac{k}{j}\right)^{q}\, \tilde\Theta_j\, \tilde\Theta_k\, \delta_{jk} + b_1 \sum_k \sum_j 
 \left(\frac{k}{j}\right)^{q}\, \tilde\Theta_j\, \tilde\Theta_k \int_{I_z}z\, \psi_k\psi_j\, d\pi(z) \\[4pt]
&\displaystyle\qquad\quad = b_0 \sum_k\sum_j \left(\frac{k}{j}\right)^{q}\, \tilde\Theta_j\, \tilde\Theta_k\, \delta_{jk}
+ b_1 \sum_k \sum_{j=k-1, k, k+1} \left(\frac{k}{j}\right)^{q}\, \tilde\Theta_j\, \tilde\Theta_k \int_{I_z}z\, \psi_k\psi_j\, d\pi(z) \\[4pt]
&\displaystyle\qquad\quad = b_0 \sum_k \tilde\Theta_k^2 +  \text{Term B}, 
\end{align*}
then
\begin{align*}
&\displaystyle |\text{Term B}| \leq
|b_1| \sum_{k=2}^K \left|\tilde\Theta_k\, \tilde\Theta_{k-1}\left(\frac{k}{k-1}\right)^{q} \int_{I_z}z\, \psi_k\psi_{k-1}\, d\pi(z)\right| \\[4pt]
&\displaystyle\qquad\qquad\quad + |b_1| \sum_{k=1}^{K-1}\left|\tilde\Theta_k\, \tilde\Theta_{k+1}\left(\frac{k}{k+1}\right)^{q} \int_{I_z}z\, \psi_k \psi_{k+1}\, d\pi(z)\right| + |b_1| \sum_{k=1}^{K} \left|\tilde\Theta_k^2 \int_{I_z}z\, \psi_k^2\, d\pi(z)\right| \\[4pt]
&\displaystyle\qquad\qquad \leq 2^q\, |b_1| \sum_{k=2}^{K} \left|\tilde\Theta_k\, \tilde\Theta_{k-1}\right| 
\left|\int_{I_z} z\, \psi_k \psi_{k-1}\, d\pi(z)\right|   \\[4pt]
&\displaystyle\qquad\qquad\quad + |b_1| \sum_{k=1}^{K-1}\left|\tilde\Theta_k\, \tilde\Theta_{k+1}\right| \left|\int_{I_z}z\, \psi_k \psi_{k+1}\, d\pi(z)\right| 
+ |b_1| \sum_{k=1}^{K}\tilde\Theta_k^2 \left|\int_{I_z} z\, \psi_k^2\, d\pi(z)\right|  \\[4pt]
&\displaystyle\qquad\qquad \leq 2^q\, |b_1|\, C_z \sum_{k=2}^{K} \left|\tilde\Theta_k\tilde\Theta_{k-1}\right|
+ |b_1|\, C_z \sum_{k=1}^{K-1} \left|\tilde\Theta_k\tilde\Theta_{k+1}\right|
+ |b_1|\, C_z \sum_{k=1}^{K}\tilde\Theta_k^2 \\[4pt]
&\displaystyle\qquad\qquad \leq 2^q\, |b_1|\, C_z\, \frac{1}{2}\left(\sum_{k=2}^K \tilde\Theta_k^2 
+\sum_{k=2}^K \tilde\Theta_{k-1}^2\right) +  |b_1|\, C_z\, \frac{1}{2}\left(\sum_{k=1}^{K-1}\tilde\Theta_k^2 + \sum_{k=1}^{K-1}\tilde\Theta_{k+1}^2\right)
+ |b_1|\, C_z \sum_{k=1}^{K}\tilde\Theta_k^2 \\[4pt]
&\displaystyle\qquad\qquad \leq (2^q+2) |b_1|\, C_z \sum_{k=1}^K \tilde\Theta_k^2, 
\end{align*}
where in the second inequality we used $\frac{k}{k-1}\leq 2$ and $\frac{k}{k+1}<1$; in the third inequality the Cauchy-Schwarz inequality is used: 
\begin{align*}
&\displaystyle \left| \int_{I_z}z\, \psi_k \psi_{k-1}\, d\pi(z)\right| \leq ||z||_{L^{\infty}} \int_{I_z}\left|\psi_k \psi_{k-1}\right| d\pi(z) \\[2pt]
&\displaystyle\qquad\qquad\qquad\qquad\quad \leq C_z \left(\int_{I_z}\psi_k^2\, d\pi(z)\right)^{1/2}
\left(\int_{I_z}\psi_{k-1}^2\, d\pi(z)\right)^{1/2} = C_z, 
\end{align*}
and $|xy|\leq\frac{1}{2}\left(|x|^2 + |y|^2\right)$ is used in the fourth inequality. 

Under the assumption on $b$ given in (\ref{BK2}), then 
\begin{align}
&\displaystyle\text{Term A} \geq b_0 \sum_{k=1}^K \tilde\Theta_k^2 - (2^q+2) |b_1|\, C_z \sum_{k=1}^K \tilde\Theta_k^2  \\[2pt]
&\displaystyle\label{TermA}\qquad\quad = \left(b_0 - (2^q+2)|b_1|\, C_z\right)\sum_{k=1}^K \tilde\Theta_k^2 \geq D(\cos\theta)\sum_{k=1}^K \tilde\Theta_k^2. 
\end{align}

Since $\mathcal M>0$, $\mathcal M^{\ast}>0$, $\phi(|v-v^{\ast}|)\geq 0$, by using (\ref{TermA}), one has Term I controlled by: 
\begin{align}
&\displaystyle \text{Term I} = -\frac{1}{4}\int\int\mathcal M\mathcal M^{\ast}\sum_k \sum_j \left(\frac{k}{j}\right)^{q}\, S_{kj}\, \tilde\Theta_j\, \tilde\Theta_k\, dv^{\ast}d\sigma dv \notag \\[4pt]
&\displaystyle\qquad\quad =  -\frac{1}{4}\int\int\mathcal M\mathcal M^{\ast}\,  \phi(|v-v^{\ast}|)\underbrace{\sum_k \sum_j \left(\frac{k}{j}\right)^{q}\, \tilde S_{kj}\, \tilde\Theta_j\, \tilde\Theta_k}_{\text{Term A}} dv^{\ast}d\sigma dv \notag \\[4pt]
&\displaystyle\qquad\quad \leq -\frac{1}{4}\sum_k \int\int\mathcal M\mathcal M^{\ast} \phi(|v-v^{\ast}|) D(\cos\theta)\, \tilde\Theta_k^2\, dv^{\ast}d\sigma dv\notag\\[4pt]
&\displaystyle\qquad\quad = -\frac{1}{4}\sum_k\, k^{2q}\int\int\mathcal M\mathcal M^{\ast} \phi(|v-v^{\ast}|) D(\cos\theta)\, \Theta_k^2\,  dv^{\ast}d\sigma dv\notag\\[4pt]
&\displaystyle\qquad\quad = \sum_k\, k^{2q}\int\int\mathcal M\mathcal M^{\ast} \underbrace{\phi(|v-v^{\ast}|) D(\cos\theta)}_{\tilde D} \Theta[\tilde h_k]\, \tilde h_k\, dv^{\ast}d\sigma dv\notag\\[4pt]
&\displaystyle\label{TermI}\qquad\quad = \sum_k\, k^{2q}\, \langle\mathcal L^{\tilde D}(h_k), \, h_k\rangle_{L^2_v}, 
\end{align}
where the equivalence between (\ref{E6a}) and (\ref{E1}) is used again in the second last row, and (\ref{E0}) is used in the last equality, 
except that now we have $\tilde D$ instead of $S_{kj}$ in the integral, with $\mathcal L^{\tilde D}$ defined as replacing $B$ by 
$\tilde D=\phi(|v-v^{\ast}|)D(\cos\theta)$ in (\ref{LB1}), with $D(\cos\theta)$
sharing the same assumption as $b(\cos\theta)$ in \cite[assumption (1.3)]{MC}, so that the coercivity property (\ref{coercivity}) still holds. 
We remark that condition (\ref{BK2}) is required for technical reasons in equations (\ref{TermA}) and (\ref{TermI}). 
It gives a restriction on the relation between $b_0$ and $b_1$ in the 
formulation (\ref{B-LN}), which is reasonable since the random perturbation part is usually not expected to be too large. 

Now integrating on $x$, one finally has
$$\text{Term I} \leq \sum_{k=1}^{K}\, k^{2q}\, \langle\mathcal L^{\tilde D}(h_k), \, h_k \rangle_{L^2} \leq 
- C_{\lambda} \sum_{k=1}^K\, ||k^q h_k^{\perp}||_{\Lambda}^2\,, $$
where $C_{\lambda}$ is a constant independent of $K$ and $z$.
\\[2pt]

The estimate on the nonlinear term $\mathcal F_k$ ($k=1, \cdots, n$) will be the same as in \cite{LiuJin-UQ}. The reason is that the upper bound for the triple index coefficient matrix that we need stays the same as before (see the inequality (5.9) in \cite{LiuJin-UQ}):
\begin{align*}
&\displaystyle |S_{mnk}| \leq (C_b + |b_1|\, C_z)\, ||\psi_n||_{L^{\infty}(z)}\, \langle |\psi_m|, \, |\psi_k|\rangle_{L^2(\pi(z))} \\[2pt]
&\displaystyle \qquad\quad \leq  (C_b + |b_1|\, C_z)\, ||\psi_n||_{L^{\infty}(z)}\, ||\psi_m||_{L^2(\pi(z))}\, ||\psi_k||_{L^2(\pi(z))} = \tilde C\, n^p\,,
\end{align*}
where $\tilde C=\mathcal O(1)$ no matter whether $|b_1|\leq \mathcal O(\varepsilon)$ or $|b_1| = \mathcal O(1)$, since $C_b=\mathcal O(1)$ will dominate
$\tilde C$.
Another difference from \cite{LiuJin-UQ} is that we no longer need to ``absorb"  the non-diagonal part of the linearized term into the non-linear term, which has the same order of coefficients ($\mathcal O(\frac{1}{\varepsilon})$ in the incompressible Navier-Stokes scaling).
The rest of the proof is the same as that for Theorem 5.1 in \cite{LiuJin-UQ}. We omit the details here.
\\[6pt]
\indent The following Theorem has the same conclusion as Theorem \ref{thm:5.1} (i.e., Theorem 5.1 in \cite{LiuJin-UQ}), namely,
the gPC based Stochastic Galerkin method for the Boltzmann equation with random inputs and both scalings of $\alpha$ is of spectral accuracy,
and the total gPC error decays exponentially in time.
The {\it significant difference} here is that we are now able to remove the ``bad"
assumption on the small $\mathcal O(\varepsilon)$ random perturbation of the collision kernel shown in (\ref{Assump-B}).
We consider the one-dimensional random space. There is no new technical difficulty on
extending the analysis to a higher dimensional case. One can refer to Remark 2.8 in \cite{ShuJin} for a discussion. 

\begin{theorem}
\label{main:thm}
We first give a summary of the assumptions needed on the collision kernel $B$: 
\begin{align}
&\displaystyle B(|v-v^{\ast}|, \cos\theta, z) = \phi(|v-v^{\ast}|)b(\cos\theta,z), \, 
\Phi(|v-v^*|)=C|v-v^*|^{\gamma}, \, \gamma\in[0,1], \, C>0, \notag\\[4pt]
&\displaystyle \forall\eta\in[-1,1], \,  |b(\eta,z)|\leq C_b, \, |\partial_{\eta} b(\eta,z)|\leq \tilde C_b, \,
|\partial_z^k b(\eta, z)|\leq C_b^{\ast}, \qquad \forall\, 0\leq k\leq r, \notag\\[4pt]
&\displaystyle b(\cos\theta, z)= b_0(\cos\theta)+  b_1(\cos\theta)z,  \qquad  |z| \leq C_z,  \notag\\[4pt]
&\displaystyle\label{b0-b1} b_0 (\cos\theta) \geq (2^q+2)\, |b_1(\cos\theta)|\, C_z + D(\cos\theta),  \\[4pt]
&\displaystyle  \inf_{\sigma_1, \sigma_2\in\mathbb S^2}\int_{\sigma_3\in\mathbb S^2}\min\{b(\sigma_1\cdot\sigma_3, z), b(\sigma_2\cdot\sigma_3, z)\}
d\sigma_3 >0 \qquad  \text{for all }z, \notag\\[4pt]
&\displaystyle \inf_{\sigma_1, \sigma_2\in\mathbb S^2}\int_{\sigma_3\in\mathbb S^2}\min\{D(\sigma_1\cdot\sigma_3), D(\sigma_2\cdot\sigma_3)\} 
d\sigma_3 >0, \notag
\end{align}
and instead of (\ref{Assump-B}) in \cite{LiuJin-UQ}, assume that 
\begin{equation} |\partial_z b| \leq R, \end{equation}
for some constant $R>0$. Let the gPC polynomial basis satisfy 
$$||\psi_k||_{L^{\infty}} \leq C k^p, \qquad 1\leq k\leq K, $$
with a parameter $p \geq 0$, and $q$ in (\ref{b0-b1}) satisfies $q>p+2$. 
The initial conditions for energy $E^K$ (defined in (\ref{Def-E})) and solution $h$ shown in (\ref{E-IC})--(\ref{h-IC}) are also satisfied. Then we have
\begin{align*}
&\displaystyle E^{K}(t)\leq \eta\, e^{-\varepsilon^{1-\alpha}\,\tau t}\,,  \\[6pt]
&\displaystyle ||h^K||_{H^s_{x,v}L_z^{\infty}} \leq \tilde\eta\, e^{-\varepsilon^{1-\alpha}\,\tau t}, \qquad ||h^K||_{H^s_{x,v}L_z^2} \leq \tilde\eta\,
e^{-\varepsilon^{1-\alpha}\,\tau t}\,,  \\[6pt]
&\displaystyle ||h-h^K||_{H_{x,v}^{s} L_z^2} \leq C\, \frac{e^{-\varepsilon^{1-\alpha}\, \tilde\tau t}}{K^r}\,.
\end{align*}
Here $\eta$, $\tilde\eta$, $\tau$, $\tilde\tau$ and $C$ are all positive constants and independent of $K$ and $z$.
\end{theorem}

\noindent {\bf Remark.} The idea of using change of variables and writing Term I in the form of (\ref{E6a}) is inspired by \cite{HY-MB},
where an explicit spectral-gap estimate and a convergence to equilibrium of the linearized multi-species Boltzmann equation was studied.
In \cite{HY-MB}, they consider a system of Boltzmann equations that models
the evolution of a dilute ideal gas composed of $n\geq 2$ different species,
\begin{align}
\label{F-MB}
\left\{
\begin{array}{l}
\displaystyle\partial_t F_i + v\cdot\nabla_x F_i = \mathcal Q_i(F), \qquad 1\leq i \leq n, \\[4pt]
\displaystyle F_i(x,v,0)=F_{I,i}(x,v), \qquad (x,v)\in {\mathbb T^3} \times {\mathbb R^3}.
\end{array}\right.
\end{align}
$\mathcal Q_i$ is the $i$-th component of the nonlinear collision operator, defined by
\begin{align*}
&\displaystyle \mathcal Q_i(f) = \sum_{j=1}^n \mathcal Q_{ij}(f_i, f_j), \\[2pt]
&\displaystyle  \mathcal Q_{ij}(f_i, f_j) =\int_{\mathbb R^3 \times \mathbb S^2} B(|v-v^{\ast}|, \cos\theta, z)
(f_i^{\prime}f_j^{\prime\ast} - f_i f_j^{\ast}) dv^{\ast}d\sigma.
\end{align*}
The linearized equation of (\ref{F-MB}) is given by
\begin{align*}
\left\{
\begin{array}{l}
\displaystyle \partial_t f_i + v\cdot\nabla_x f_i = L_i(f), \\[4pt]
\displaystyle f_i(x,v,0)=f_{I,i}(x,v), \qquad (x,v)\in {\mathbb T^3} \times {\mathbb R^3},
\end{array}\right.
\end{align*}
where $F_i = \mathcal M_i + \sqrt{\mathcal M_i}\, f_i$, and
\begin{align*}
&\displaystyle L_i(f) = \sum_{j=1}^n L_{ij}(f_i, f_j), \\[2pt]
&\displaystyle  L_{ij}(f_i, f_j) = \int_{\mathbb R^3 \times \mathbb S^2} B_{ij}\, \mathcal M_i^{1/2}\mathcal M_j^{\ast}
\left(h_i^{\prime} + h_j^{\prime, \ast} - h_i - h_j^{\ast}\right) dv^{\ast}d\sigma, \qquad h_i:=\frac{f_i}{\sqrt{\mathcal M_i}}\,. 
\end{align*}

One important step relevant to our analysis is the following idea. Split the operator $L=L^m + L^b$ with $L^m=(L_1^m, \cdots, L_n^m)$ and  $L^b=(L_1^b, \cdots L_n^b)$, given by
$$ L_i^m(f_i) = L_{ii}(f_i, f_i), \qquad L_i^b = \sum_{j \neq i} L_{ij}(f_i, f_j), $$
then $(f, L^m(f))_{L^2_v}$ and $(f, L^b(f))_{L^2_v}$ can be written in bilinear forms with squares:
\begin{align*}
&\displaystyle (f, L^m(f))_{L^2_v} = - \frac{1}{4} \sum_{i=1}^n \int_{\mathbb R^6\times \mathbb S^2}\,  B_{ii}\Delta_i[h_i]^2\, 
\sqrt{\mathcal M_i} \sqrt{\mathcal M_i^{\ast}}\, dv dv^{\ast}d\sigma,
\\[4pt]
&\displaystyle (f, L^b(f))_{L^2_v} = - \frac{1}{4}\sum_{i=1}^n \sum_{j\neq i} \int_{\mathbb R^6\times \mathbb S^2}\,  
B_{ij} A_{ij}[h_i, h_j]^2\, \sqrt{\mathcal M_i}\sqrt{\mathcal M_j^{\ast}}\, dv dv^{\ast}d\sigma,
\end{align*}
where $\Delta_i[h_i] := h_i^{\prime} + h_i^{\prime\ast} - h_i - h_i^{\ast}$,
$A_{ij}[h_i, h_j]:= h_i^{\prime} + h_j^{\prime\ast} - h_i - h_j^{\ast}$.
Since $B_{ij}\geq 0$ for all $1\leq i, j \leq n$, thus $(f, L^m(f))_{L^2_v}\leq 0$ and $(f, L^b(f))_{L^2_v}\leq 0$.

\bibliographystyle{siam}
\bibliography{Boltzmann_gPC.bib}

\medskip
Received xxxx 20xx; revised xxxx 20xx.
\medskip

\end{document}